# Ax, 3 polyominoes for tiling the plane non-periodically

Vincent Van Dongen[1] and Pierre Gradit[2]

[1] Independent researcher, Canada; vincent.vandongen@gmail.com
[2] Independent researcher, France; pierre.gradit@gmail.com

**Abstract**

How do people come up with new sets of tiles including new tile shapes that would only tile non-periodically? This paper presents our graphical journey in tilings and provides a new set of three polyominoes named Ax for its relationship with Ammann A4.

## Introduction

There are numerous aperiodic set of tiles in the literature [1] [2]. We wonder how these tilesets can be created. This paper unveils a graphical process by which we derive a new set of 3 tiles from the *"Fourth Ammann's Set of Tiles"* (A4) [2]. Even if our final goal is to provide an aperiodic set, our exploration method is based on a tileset where each tile embeds a representation of its supertile. This kind of representation was already used in [2] for A4. We call it a *"self-ruling tile"*. See Figure 1(b). Note that for A4, no supertile is associated to the key tile. Using self-ruling tiles A and B for A4, one can create larger tiles iteratively. Also note that both sides of Tiles A and B are used in the self-ruling tiles.

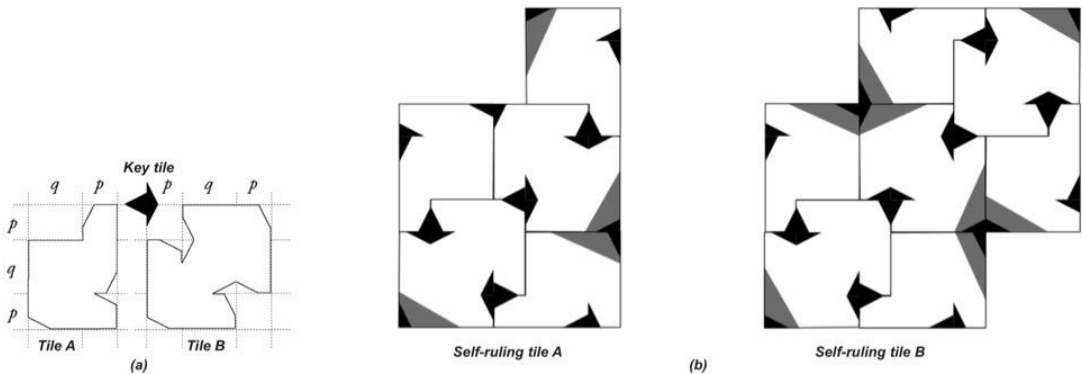

**Figure 1:** (a) *Ammann A4 consists of tile A, tile B and key tile. (b) The self-ruling tiles of A4.*

In the next section, we will present our starting point, a set of self-ruling polyominoes whose boundaries are the same as Ammann A4 with p=3 and q=2. As a reminder, a polyomino is a geometric figure formed by unit squares joined edge to edge.

## Our starting point: a set of self-ruling polyominoes

We were exploring the use of Hitomezashi marking on Ammann A4 tile. This type of marking was already used with success on Ammann-Beenker tiling [3].

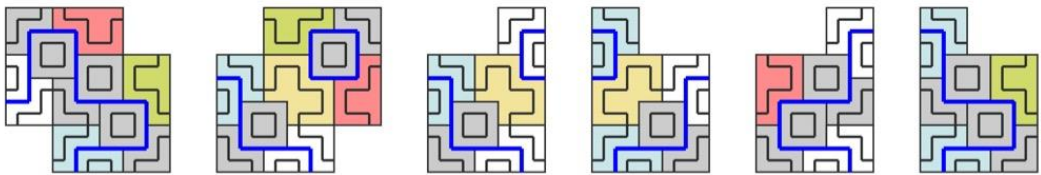

**Figure 2:** *A set of self-ruling polyominoes with Hitomezashi marking.*

After many trials and errors of marking tiles A and B of A4, we came up with the set of self-ruling polyominoes of Figure 2. The Hitomezashi marks are the stitches of unit size in black on the tiles.

Is there a way to reduce the number of self-ruling tiles? In the next section, we will provide the steps to get there. But before, we first explain how to read these particular self-ruling tiles. As an example, Figure 3(c) is one of the polyominoes of Figure 2. Figure 3(b) highlights its blue curves. Figure 3(a) shows the tile that represents it (with the same curves but in black). In Figure 3(d), the supertile is exploded into tiles for its expansion. Figure 3(e) shows the replacement of the tiles with the corresponding supertiles. Figure 3(f) shows the result of that replacement. This is an iterative process to be done over and over, on all tiles, to create larger and larger tilings.

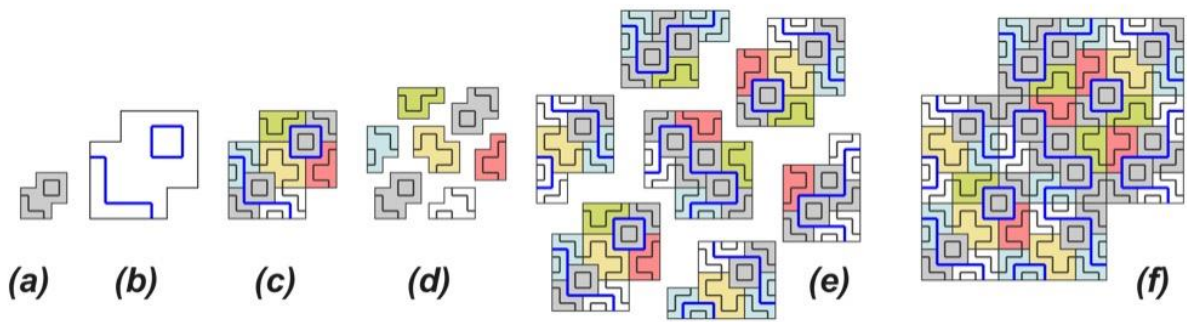

**Figure 3:** *Illustration of the use of a self-ruling polyomino.*

It is worth saying a few words on the blue curves that are on the supertiles. They are similar to Hitomezashi threads except that the stitches are not of length 1 but rather of length that is either 2 or 3. In fact, we observed that they are on a grid similar to Ammann bar grid although its orientation is different than the Ammann bar grids of A4 found in the literature [4]. At the end of the paper, more will be said on this.

## Reducing the number of self-ruling polyominoes

We will now show how we were able to reduce the set of self-ruling tiles of Figure 2. Looking at many large tilings (not shown here due to lack of space), we first conjecture that the white tile and the blue one are always grouped together (these groups appear in Figure 3(f) for example. As a result, we can replace both tiles by a new one. And this lead us to the new set of tiles.

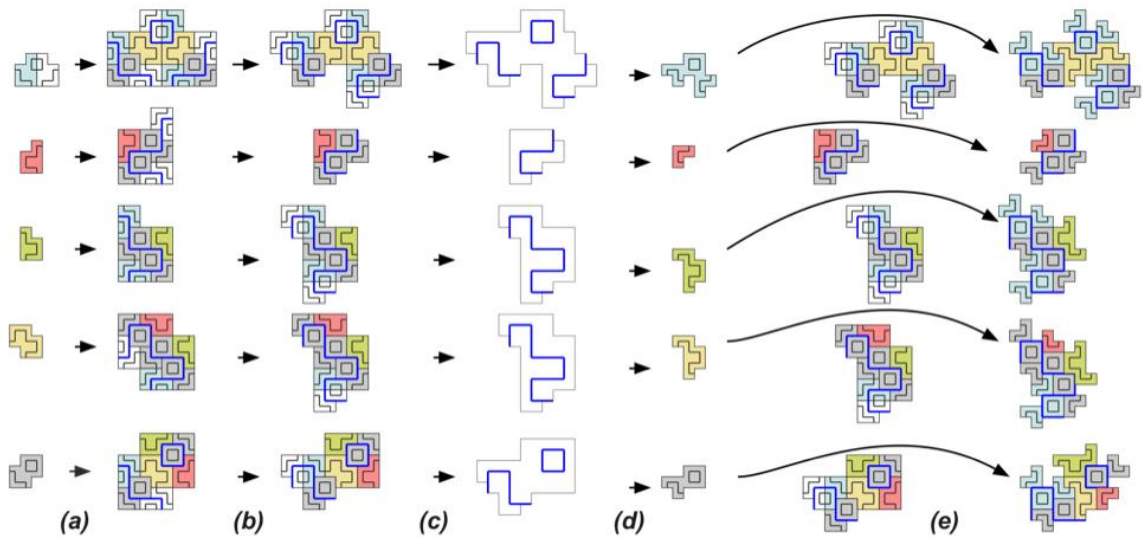

**Figure 4:** *The grouping of 2 tiles that lead to a temporary set of new self-ruling polyominoes.*

See at the top of Figure 4(a) the grouping of the blue and white tile and its associated supertile. At Figure 4(b), this group is enforced everywhere in the tile including at its borders. Figure 4(c) shows the blue curves of the new supertiles and Figure 4(d) shows the associated new tiles. Figure 4(e) shows the last step that consists of rebuilding the new self-ruling tiles based on these new tiles.

By observing carefully the new tileset, we can now optimize it. Clearly, two tiles are identical, the ones on 3$^{rd}$ and 4$^{th}$ line of Figure 4(e). Also, the first line of 4(e) is no longer necessary as it can created by means of 2 other tiles, the ones of line 2 and line 5. This further optimization is shown in Figure 5. Figure 5(a) shows this new rewriting leading to the new self-ruling tileset of Figure 5(b).

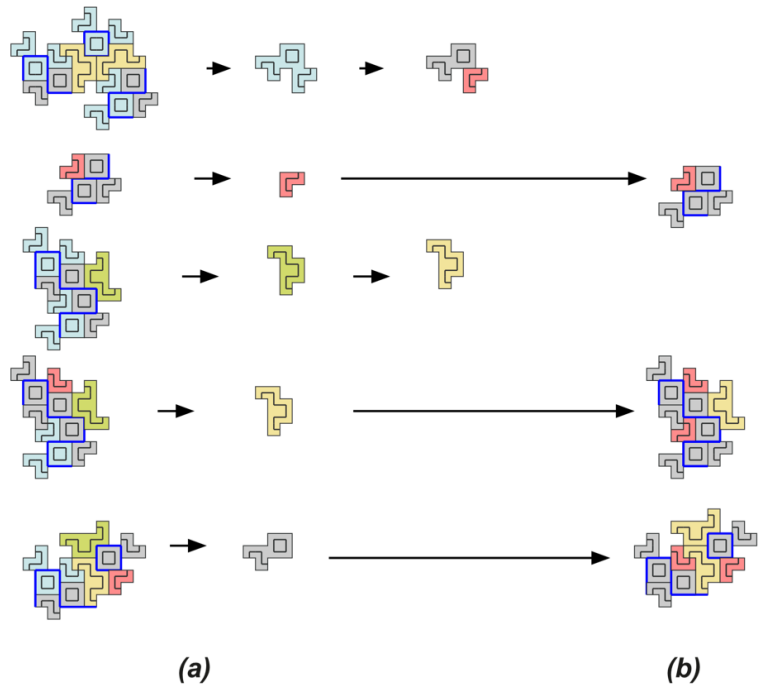

**Figure 5:** *Tileset optimization: going from (a) 5 to (b) 3 self-ruling polyominoes.*

## Testing the new tileset Ax

We give the name Ax to our new tileset of 3 self-ruling polyominoes. Figure 6 shows the initial test of Ax.

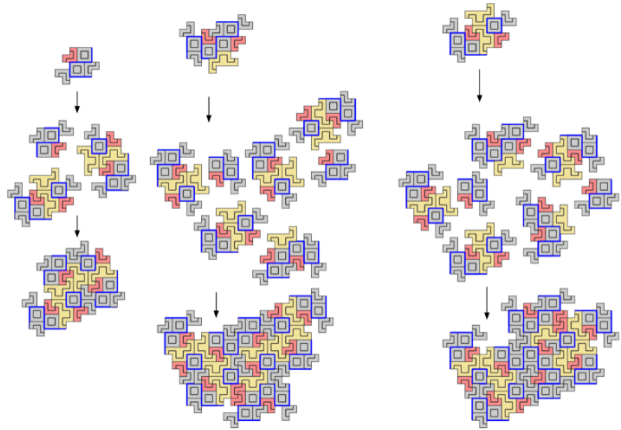

**Figure 6:** *The 3 self-ruling polyominoes of Ax and their expansion after one iteration.*

Here we want to make sure that no overlap nor holes are created. Clearly, it works at this iteration. We then repeat the process to create tilings. After 4 or 5 iterations, if it's still ok, we can conjecture that it works. See Figure 7 for an example of tiling created. We observe that closed curves are Fibonacci snowflakes.

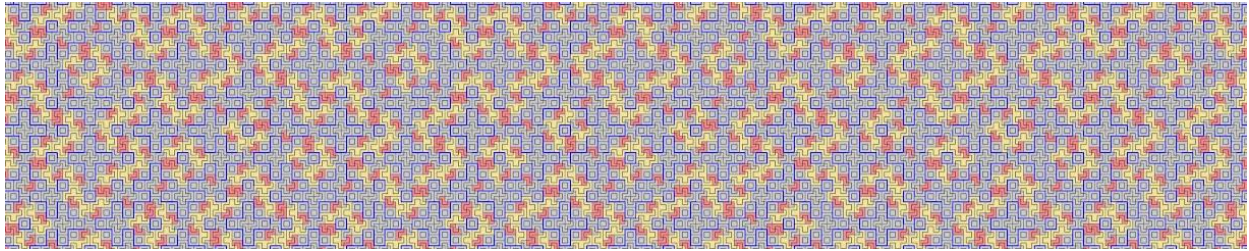

**Figure 7:** *Example of tiling Ax (cropped image).*

### About the potential aperiodicity of tileset Ax

Proving that Ax is aperiodic is outside the scope of this paper. Here, we can give some initial hints. First, it is known that A3, another Ammann's tiling, is aperiodic, Ammann bars enforcing the aperiodicity of its tiles [1]. To recall, Ammann A3 is an aperiodic set of 3 tiles with Ammann bars at 45 degrees. In the same way, we can hope the tileset of Figure 2 to be aperiodic. Figure 8 shows the way the Ammann bars can be enforced with tile shape only, using one of its supertiles as an example.

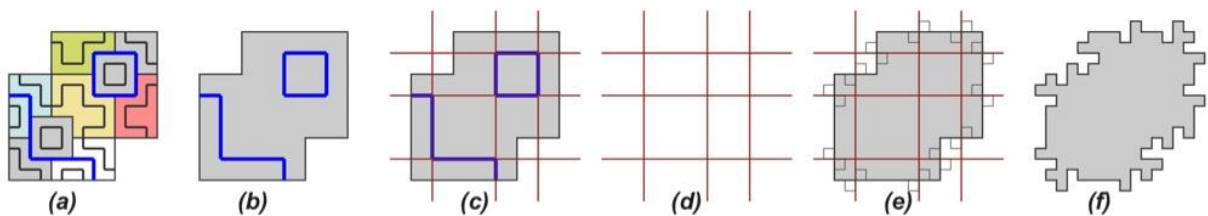

**Figure 8:** *Example of implementation of the Ammann bars by shape only.*

In the same manner, Ammann bars are on Ax and can be enforced by shape only. Based on this, there is a chance that Ax is aperiodic. But a formal proof based on these observations is yet to come.

### Conclusion

Finding a new tileset is an interesting journey. We presented here our experience in this field. As shown in the paper, it is primarily about observing patterns and rewriting tiles into new ones. Our journey for this new tileset, named Ax, is not over. We are currently developing tools to prove that Ax is aperiodic. Also we already found many other sets of 3 polyominoes that can be considered as variations of Ax. The one shown here is made of self-ruling tiles and the patterns created are Fibonacci snowflakes.